\documentclass[a4paper,12pt]{article}
\usepackage[cp1251]{inputenc}
\usepackage[english]{babel}
\usepackage{amsfonts,amssymb,amsmath}
\usepackage[dvips]{graphics}
 \newcommand{\re}{\mathop{\rm Re}\nolimits}

\begin{document}

Asymptotic solutions of the boundary value problems for the singularly perturbed differential algebraic equations with a turning point

\vspace{5mm}

Petro Samusenko\\ Department of Mathematical Analysis and Probability Theory,
Faculty of Physics and Mathematics, National Technical University of Ukraine
“Igor Sikorsky Kyiv Polytechnic Institute”, Kyiv, Ukraine\\ psamusenko@ukr.net

\vspace{5mm}

This paper deals with the boundary value problems for the singularly pertur\-bed differential-algebraic system of equations. The case of turning points has been studied. The sufficient conditions for existence and uniqueness of the solution of the boundary value problems for DAEs have been found. The technique of constructing the asymptotic solutions has been developed.

\vspace{5mm}


\vspace{5mm}


\vspace{5mm}

\section{Introduction}

There are a lot of various processes in practice which are described by a differential-algebraic system of equations (DAEs)
\begin{equation}A(t)\frac{dx}{dt}=f(x,t),\;t\in[0;T],\label{equ:1}\end{equation}
where $x=x(t)$ is an $n$-dimensional vector function, and $A(t)$ is a square matrix of order $n$. Thus, the problems of the theory of electric circuits  \cite{cho, chu, dem, gun, ria1, tish}, theory of optimal control \cite{bal, camp1, ger, gri, kun}, neural networks \cite{camp2, new, ria2}, hydrodynamics \cite{camp3, camp4, camp5}, robotics \cite{camp6, new, spi1, spi2}, chemistry \cite{kum} lead to the need of studying certain properties of the solutions of linear DAEs
\begin{equation}A(t)\frac{dx}{dt}=B(t)x+f(t),\;t\in[0;T].\label{equ:2}\end{equation}

The DAEs (\ref{equ:2}) began to be intensively studied since the late 1970's, although the first results on DAEs were obtained much earlier.
In particular, Luzin \cite{luz} and Gantmacher \cite{gan} have found necessary and sufficient conditions for sol\-va\-bi\-lity of DAEs with constant coefficients, and they have proposed some approaches of constructing their solutions. It should be noted that Gantmacher's algorithm for constructing of the particular solutions of DAEs
\begin{equation}A\frac{dx}{dt}=Bx+f(t)\label{equ:3}\end{equation}
was based on the idea of reduction of a pencil $A-\lambda B$ to Kronecker normal form. At the same time similar technique could not be used to DAEs with variable coefficients
\begin{equation}\label{equ:4}
A(t)\frac{dx}{dt}=B(t)x+f(t)\end{equation} in general case, since
its application can cause changing the Kronecker's form of
a pencil $A(t)-\lambda B(t)$.

A more general approach to the study of DAEs is given in \cite{camp7, camp8, camp9}, where the concept of the central canonical form was introduced.
This approach can be also used to the study of systems with variable coefficients.

If by means of nonsingular transformations the system (\ref{equ:4}) can be written in the form
\begin{equation}\begin{pmatrix}
  I_{n-s} & 0 \\
  0 & N_s(t)
\end{pmatrix}\frac{dy}{dt}=\begin{pmatrix}
  M(t) & 0 \\
  0 & I_s
\end{pmatrix}y+h(t),\label{equ:5}\end{equation}
where $I_s$ and $I_{n-s}$ are identity matrices of orders $s$ and $n-s$, respectively, and
$N_s(t)$ is a nilpotent lower (or upper) triangular matrix
\cite{camp9}, then the system (\ref{equ:5}) is called the standard
canonical form of system (\ref{equ:4}). Note that, when matrix $N_s(t)$ is a matrix with constant coefficients, the system (\ref{equ:5}) is called the strong standard
canonical form \cite{gea,pet}.

The effective methods of transforming DAEs to a standard canonical form are presented by
Boyarintsev \cite{boy}, Campbell \cite{camp9}, Petzold \cite{pet}, Samoilenko, Shkil' and Yakovets
\cite{samo} if the rank of the matrix $A(t)$ is constant on the interval $[0;T]$. It allows to find the general solution of the system
(\ref{equ:4}) and to study then Cauchy problem, boundary value
problems, and others \cite{camp9,camp8,boy}.

Another method to deal with DAEs is based on the concept of the tractability index DAEs \cite{grie, marz1, marz2}. In this case, using the projectors technique, it is possible to find the solution of DAEs in the form of the sum of several terms, one of which is the solution of inherent regular ODEs, and the other terms are taken from the corresponding algebraic systems \cite{marz3, ria3, lam1}. The projectors technique for finding solutions of DAEs with critical points is generalized in \cite{marz4, marz5, koch}.

By means of the linearization procedure, the methods for solving linear DAEs described above can be used for solving of nonlinear DAEs \cite{marz6, lam1}.

Using a concept of the tractability index, the numerical methods for solving differential-algebraic system have been developed in the papers by Gear and Petzold \cite{gea}, Griepentrog and M\"arz \cite{grie}, Brenan, Campbell and Petzold \cite{bre}.

Now there are several approaches to solving boundary value problems for DAEs. Thus, sufficient conditions for the existence and uniqueness of two-point boundary value problems for DAEs index 1 and 2 have been found in \cite{marz7, marz8, len, lam2}. Using the projectors technique, analytical methods for finding solutions of these boundary value problems are developed. Similar approaches to solving multipoint boundary value problems are given in \cite{marz7, anh}. The standard canonical form is used to solving a multipoint boundary value problems in \cite{boi1, boi2}.

One of the efficient methods of integration of DAEs is the
perturbation method \cite{camp7, marz8,han1,han2,boy} according to which instead of system
\begin{equation} A(t)\frac{dx}{dt}=B(t)x+f(t),\;t\in[0;T],\label{equ:6}\end{equation}
one can consider a perturbed system \begin{equation} (A(t)+\varepsilon A_1(t))
\frac{dx}{dt}=B(t)x+f(t),\label{equ:7}\end{equation}
where $\varepsilon$ is a small positive parameter.
The matrix $A_1(t)$ is chosen so that $det(A(t)+\varepsilon A_1(t))\not=0$ for all $t\in[0;T]$, and all sufficiently small
$\varepsilon >0$. Thus, the system (\ref{equ:7}) is a regular DAEs, known methods can be used to solve it.

The conditions under which the solutions of the system (\ref{equ:7}) converge to the
corresponding solutions of system (\ref{equ:6}), as
$\varepsilon\rightarrow 0+$ can be found in \cite{marz8,han1,han2,boy}.

Consider the more general system than (\ref{equ:7}), namely,
\begin{equation}\varepsilon
A(t,\varepsilon)\frac{dx}{dt}=B(t,\varepsilon)x+f(t,\varepsilon),\;t\in[0;T],\label{equ:8}\end{equation}
where $A(t,\varepsilon)$, $B(t,\varepsilon)$ are square matrices of orders $n$, $f(t,\varepsilon)$ is an $n$-dimensional vector
possessing uniform asymptotic expansions of the following form
\[A(t,\varepsilon)=\sum_{k=0}^{\infty}\varepsilon^k A_k(t),\; B(t,\varepsilon)=\sum_{k=0}^{\infty}\varepsilon^k B_k(t),\; f(t,\varepsilon)=\sum_{k=0}^{\infty}\varepsilon^k f_k(t)\] with real or complex-valued infinitely differentiable coefficients.
The solutions of such singular perturbed DAEs have a number of specific features in comparison with the solutions of system (\ref{equ:6}).
Samoilenko, Shkil' and Yakovets have shown that under certain conditions for perturbed matrices a homo\-ge\-neous system
\begin{equation}\varepsilon
A(t,\varepsilon)\frac{dx}{dt}=B(t,\varepsilon)x\label{equ:9}\end{equation}
has two types of linear independent formal solutions corresponding
to finite or infinite elementary divisors of a pencil
$B_0(t)-\lambda A_0(t)$ \cite{samo}. Moreover, their linear
combination is the formal general solution of system (\ref{equ:9}). It
should be noted that in case of multiple elementary divisors of a pencil
$B_0(t)-\lambda A_0(t)$ the asymptotic expansions of the solutions of
the system (\ref{equ:9}) can be constructed in some fractional powers
of small parameter $\varepsilon$, where values of powers of
$\varepsilon$ depends on the multiplicity of the roots of a
characteristic equation
$$det(B_0(t)-\lambda A_0(t))=0$$ with corresponding elementary divisors, as well
as on perturbed coefficients of the system (\ref{equ:9}).

Note that the transformation of the system (\ref{equ:8}) using the projectors technique is not rational, because such the transformation significantly changes the characteristics of the system, for example, the structure of elementary divisors of the pencil $B_0(t)-\lambda A_0(t)$.

Using the results of asymptotic analysis of singularly perturbed DAEs (\ref{equ:9}), the two-point boundary value problem
\begin{equation} \varepsilon A(t,\varepsilon)\frac{dx}{dt}=f(x,t,\varepsilon),\;t\in[0;T],\;\varepsilon\in(0;\varepsilon_0],\label{equ:10}\end{equation}
\begin{equation} Mx(0,\varepsilon)+Nx(T,\varepsilon)=d(\varepsilon),\label{equ:11}\end{equation} has been studied in \cite{vira}, where $x(t,\varepsilon)$ is an $n$-dimensional vector, $A(t,\varepsilon)$ is a square matrix of order $n$, $M$, $N$ are rectangular matrices of dimension $m\times n$, $f(x,t,\varepsilon)$ and $d(\varepsilon)$ are vectors of dimension $n$ and $m$, respectively, $\varepsilon$ is a small parameter.

Note that in the case of the linear vector function $f(x,t,\varepsilon)$ the boundary value problem (\ref{equ:10}), (\ref{equ:11}) was comprehensively studied by Yakovets and Vira in \cite{yak,vira1}. These researchers have found the conditions for the existence and uniqueness of the solution of the problem (\ref{equ:10}), (\ref{equ:11}) and have constructed its asymptotic expansion in powers of small parameter.

It should be noted that both the system (\ref{equ:8}) and (\ref{equ:10}) were studied in the absence of turning points. In addition, the constancy of the rank of the matrix $A(t,\varepsilon)$ on the set $[0;T]\times (0;\varepsilon_0]$ was assumed.

The present paper deals with the boundary value problem (\ref{equ:10}), (\ref{equ:11}) with turning point. It is also assumed that the rank of the matrix $A(t,\varepsilon)$ changes on the set $[0;T]\times (0;\varepsilon_0]$.

We find function represented as asymptotic expansion in a small para\-meter that satisfy the boundary value problem (\ref{equ:10}), (\ref{equ:11}) with certain accu\-ra\-cy. We need to define the form of this function to propose a recurrent algorithm, determining all of the terms of the corresponding asymptotic solution, and, in addition, to evaluate the accuracy with which the asymptotic approximations satisfy the boundary value problem.

The solution of the boundary value problem (\ref{equ:10}), (\ref{equ:11}) is constructed by the method of boundary functions \cite{vas}. According to this technique, the formal solution of the problem (\ref{equ:10}), (\ref{equ:11}) can be found as a sum of a regular series and two boundary layer series \cite{vas}. The presence of boundary layer series allows us to construct a uniform asymptotic solution of the problem (\ref{equ:10}), (\ref{equ:11}) on the segment $[0;T]$. The algorithm for constructing asymptotic solutions of boundary value problems developed in this paper is analogous to the algorithms for asymptotic integration of singularly perturbed problems \cite{vas}.

This  paper  is  organized  as  follows.  In  Sec.  II, an  algorithm  for  constructing  an  formal solution  to  the boundary value problem (\ref{equ:10}), (\ref{equ:11}) is  proposed  and  described  in  detail. In Sec. III, we prove the asymptotic nature of the constructed formal solution.

\section{Formal solutions}

Assume that the following conditions are satisfied:\\ 1. Elements of matrix $A(t,\varepsilon)$ have an infinite number of continuous partial derivatives with respect to $t$ and $\varepsilon$ ($A(t,\varepsilon)\in C^{\infty}(G)$) in some domain \[G=\{(t,\varepsilon):0\leq t\leq T,\,0\leq\varepsilon\leq\varepsilon_0\}\]
2. Components of vector-function $f(x,t,\varepsilon)$ have an infinite number of continuous partial derivatives with respect to $x$, $t$ and $\varepsilon$ ($f(x,t,\varepsilon)\in C^{\infty}(K)$) in some domain \[K=\{(x,t,\varepsilon):||x||<+\infty,\, 0\leq t\leq T,\,0\leq\varepsilon\leq\varepsilon_0\}.\] 3. Equation $f(x,t,0)=0$ has the solution $x=\overline{x}_0(t)$, which satisfies the conditions:\\(i) $\overline{x}_0(t)\in C[0;T]$;\\ (ii) the root $x=\overline{x}_0(t)$ is isolated on the segment $[0;T]$, that is, there is such $\eta>0$, that $f(x,t,0)\not=0$ when $0<||x-\overline{x}_0(t)||<\eta$, $t\in [0;T]$.\\4. $det A(0,0)=0$.\\5. Pencil of matrices $f_x'(\overline{x}_0(0),0,0)-\eta A(0,0)$ is regular, has two eigenvalues $\eta_1(0,0)$, $\eta_2(0,0)$, and two finite elementary divisors $(\eta-\eta_1(0,0))^p$, $(\eta-\eta_2(0,0))^q$, furthermore $p+q=n-1$. Here $f_x'(\overline{x}_0(t),t,0)=\left(\frac{\partial f_i(\overline{x}_0(t),t,0)}{\partial x_j}\right)_{i,j=\overline{1,n}}$.\\
6. $\re \eta_1(0,0)>0$, $\re \eta_2(0,0)<0$.\\
7. Pencil of matrices $f_x'(\overline{x}_0(t),t,0)-w A(t,0)$ is regular for all $t\in(0;T]$, and has $n$ distinct eigenvalues $w_i(t,0)$ such that
$w_i(t,0)\not=w_j(t,0)$, $t\in(0;T]$, $i,j=\overline{1,n}$.\\
8. $\re w_i(t,0)>0$, $t\in(0;T]$, $i=\overline{1,p+1}$, and $\re w_i(t,0)<0$, $t\in(0;T]$, $i=\overline{p+2,n}$.

The conditions (4), (5) implies that pencil of matrices $f_x'(\overline{x}_0(0),0,0)-\eta A(0,0)$
have one infinite elementary divisor of multiplicity 1.

Formal solution of the problem (\ref{equ:10}), (\ref{equ:11}) we will find in the form \begin{equation}x(t,\varepsilon)=\overline{x}(t,\varepsilon)+\Pi x(\tau,\varepsilon)+ Qx(\xi,\varepsilon),\label{equ:12}\end{equation} where $\overline{x}(t,\varepsilon)=\sum\limits_{k=0}^{\infty} \varepsilon^k\overline{x}_k(t)$ is a regular part of the asymptotics, $\Pi x(\tau,\varepsilon)=$ \linebreak[4]$=\sum\limits_{k=0}^{\infty} \varepsilon^k\Pi_k x(\tau)$, $\tau=\frac{t}{\varepsilon}$, and $Q x(\xi,\varepsilon)=\sum\limits_{k=0}^{\infty} \varepsilon^kQ_k x(\xi)$, $\xi=\frac{t-T}{\varepsilon}$, is a singular part of the asymptotics.

Substituting representation (\ref{equ:12}) in system (\ref{equ:10}), we get \[\varepsilon A(t,\varepsilon)\frac{d\overline{x}(t,\varepsilon)}{dt}+A(\varepsilon\tau,\varepsilon)\frac{d\Pi x(\tau,\varepsilon)}{d\tau}+A(\xi\varepsilon+T,\varepsilon)\frac{dQx(\xi,\varepsilon)}{d\xi}=\] \[=f(\overline{x}(t,\varepsilon)+\Pi x(\tau,\varepsilon)+Qx(\xi,\varepsilon),t,\varepsilon).\]

Then we find the functions $\overline{x}(t,\varepsilon)$, $\Pi x(\tau,\varepsilon)$, $Qx(\xi,\varepsilon)$, solving the following systems
\begin{equation}\varepsilon A(t,\varepsilon)\frac{d\overline{x}}{dt}=\overline{f}(t,\varepsilon),\label{equ:13}\end{equation}
\begin{equation}A(\varepsilon\tau,\varepsilon)\frac{d\Pi x}{d\tau}=\Pi f(\tau,\varepsilon),\label{equ:14}\end{equation}
\begin{equation}A(\xi\varepsilon+T,\varepsilon)\frac{dQx}{d\xi}=Qf(\xi,\varepsilon),\label{equ:15}\end{equation}
where \[\overline{f}(t,\varepsilon)=f(\overline{x}(t,\varepsilon),t,\varepsilon),\] \[\Pi f(\tau,\varepsilon)=f(\overline{x}(\varepsilon\tau,\varepsilon)+\Pi x(\tau,\varepsilon),\varepsilon\tau,\varepsilon)-
f(\overline{x}(\varepsilon\tau,\varepsilon),\varepsilon\tau,\varepsilon),\]
\[Q f(\xi,\varepsilon)=f(\overline{x}(\xi\varepsilon+T,\varepsilon)+Q x(\xi,\varepsilon),\xi\varepsilon+T,\varepsilon)-
f(\overline{x}(\xi\varepsilon+T,\varepsilon),\xi\varepsilon+T,\varepsilon).\]

Let \[\overline{f}(t,\varepsilon)=\sum\limits_{k=0}^{\infty}\varepsilon^k\overline{f}_k(t),\;\Pi f(\tau,\varepsilon)=\sum\limits_{k=0}^{\infty}\varepsilon^k\Pi_kf(\tau),\;Qf(\xi,\varepsilon)=\sum\limits_{k=0}^{\infty}\varepsilon^kQ_kf(\xi).\]
Here, in particular \[\overline{f}_0(t)=f(\overline{x}_0(t),t,0),\;\Pi_0f(\tau)=f(\overline{x}_0(0)+\Pi_0x(\tau),0,0)-f(\overline{x}_0(0),0,0),\]
\[Q_0f(\xi)=f(\overline{x}_0(T)+Q_0x(\xi),T,0)-f(\overline{x}_0(T),T,0),\]
\[\overline{f}_k(t)=f_x'(\overline{x}_0(t),t,0)\overline{x}_k(t)+\overline{g}_k(t),\]
\[\Pi_kf(\tau)=f_x'(\overline{x}_0(0)+\Pi_0x(\tau),0,0)\Pi_kx(\tau)+g_k(\tau),\]
\[Q_kf(\xi)=f_x'(\overline{x}_0(T)+Q_0x(\xi),T,0)Q_kx(\xi)+h_k(\xi),\;k\in N,\]
the functions $\overline{g}_k(t)$, $g_k(\tau)$ and $h_k(\xi)$  are expressed recursively through $\overline{x}_i(t)$, $\Pi_ix(\tau)$ and $Q_i(\xi)$, $i<k$.

The boundary conditions (\ref{equ:11}) can be written as
\begin{equation}
  M(\overline{x}(0,\varepsilon)+\Pi x(0,\varepsilon)+Q x(-T/\varepsilon,\varepsilon))+N(\overline{x}(T,\varepsilon)+\Pi x(T/\varepsilon,\varepsilon)+Q x(0,\varepsilon))=d(\varepsilon).
\label{equ:16}\end{equation}

The conditions (1), (2) implies that the following formal expansions hold:
\[A(t,\varepsilon)=\sum\limits_{k=0}^{\infty}\varepsilon^kA_k(t)\equiv \sum\limits_{k=0}^{\infty}\varepsilon^k\frac{1}{k!}\frac{\partial^k A(t,0)}{\partial\varepsilon^k},\]
\[A(\varepsilon\tau,\varepsilon)=\sum\limits_{k=0}^{\infty}\varepsilon^k\Pi_kA(\tau)\equiv \sum\limits_{k=0}^{\infty}\varepsilon^k\sum_{i=0}^k\frac{\tau^{k-i}}{i!(k-i)!}\frac{\partial^k A(0,0)}{\partial t^{k-i}\partial\varepsilon^i},\]
\[A(\xi\varepsilon+T,\varepsilon)=\sum\limits_{k=0}^{\infty}\varepsilon^kQ_kA(\xi)\equiv \sum\limits_{k=0}^{\infty}\varepsilon^k\sum_{i=0}^k\frac{\xi^{k-i}}{i!(k-i)!}\frac{\partial^k A(T,0)}{\partial t^{k-i}\partial\varepsilon^i},\]
\[d(\varepsilon)=\sum\limits_{k=0}^{\infty}\varepsilon^kd_k.\]

Let us equate coefficients of like powers of $\varepsilon$ in (\ref{equ:13}) -- (\ref{equ:15}). For the leading terms of the asymptotics ($\overline{x}_0(t)$, $\Pi_0 x(\tau)$ and $Q_0 x(\xi)$), we obtain
\begin{equation}f(\overline{x}_0(t),t,0)\equiv 0,\label{equ:17}\end{equation}
\begin{equation}A(0,0)\frac{d\Pi_0 x}{d\tau}=f(\overline{x}_0(0)+\Pi_0x,0,0)-f(\overline{x}_0(0),0,0),\label{equ:18}\end{equation}
\begin{equation}A(T,0)\frac{dQ_0x}{d\xi}=f(\overline{x}_0(T)+Q_0x,T,0)-f(\overline{x}_0(T),T,0).\label{equ:19}\end{equation}

In view of condition 3 systems (\ref{equ:18}), (\ref{equ:19}) will have the form
\begin{equation}A(0,0)\frac{d\Pi_0 x}{d\tau}=f(\overline{x}_0(0)+\Pi_0x,0,0),\label{equ:20}\end{equation}
\begin{equation}A(T,0)\frac{dQ_0x}{d\xi}=f(\overline{x}_0(T)+Q_0x,T,0).\label{equ:21}\end{equation}

Conditions 4, 5 guarantees the existence such nonsingular matrices $P$, $Q$, that \[PA(0,0)Q=H(0,0)\equiv diag\{0,I_{n-1}\},\] \[Pf_x'(\overline{x}_0(0),0,0)Q=\Omega(0,0)\equiv diag\{1,\Lambda_+(0,0),\Lambda_-(0,0)\},\]
where $\Lambda_+(0,0)=\eta_1(0,0)I_p+N_p$, $\Lambda_-(0,0)=\eta_2(0,0)I_q+N_q$, $I_p$ is the identity matrix of order $p$, $N_p$ is the square matrix of order $p$ such that
\[N_p=\begin{pmatrix}
  0 & 1 & 0 & \ldots & 0 \\
  0 & 0 & 1 & \ldots & 0 \\
  \hdotsfor{5}\\
  0 & 0 & 0 & \ldots & 1 \\
    0 & 0 & 0 & \ldots & 0 \\
  \end{pmatrix}\]
\cite{sib,samu}.
The matrices $I_{n-1}$, $I_q$ and $N_q$ are defined similarly. Without loss of generality we may assume that $A(0,0)=H(0,0)$ and $f_x'(\overline{x}_0(0),0,0)=\Omega(0,0)$.

Then the system (\ref{equ:20}) can be written in the form
\begin{equation}f_1(\overline{x}_0(0)+\Pi_0x,0,0)=0\label{equ:22}\end{equation}
\begin{equation}\frac{d\Pi_{02} x}{d\tau}=f_2(\overline{x}_0(0)+\Pi_0x,0,0),\label{equ:23}\end{equation}
where $\Pi_{01}x$ and $f_1(\overline{x}_0(0)+\Pi_0x,0,0)$ are the first components of the vectors $\Pi_0x$ and $f(\overline{x}_0(0)+\Pi_0x,0,0)$, respectively, and through $\Pi_{02}x$, $f_2(\overline{x}_0(0)+\Pi_0x,0,0)$ we denote the vectors containing other components of vectors $\Pi_0x$ and $f(\overline{x}_0(0)+\Pi_0x,0,0)$.

Further we make the following assumptions.
\\9. The equation $f_1(\overline{x}_0(0)+\Pi_0x,0,0)=0$ has the solution $\Pi_{01}x=\Pi_{01}x(\Pi_{02}x)$, which are continuous in the field of parameters change $\Pi_{02}x$, and \[\Pi_{01}x(\Pi_{02}x)\rightarrow 0,\;\Pi_{02}x\rightarrow 0. \]
10. The system (\ref{equ:23}) has the solution $\Pi_{02}x=\Pi_{02}x(\tau)$, such that
\begin{equation}\Pi_{02}x(\tau)\rightarrow 0,\;\tau\rightarrow+\infty,\label{equ:24}\end{equation}
and
\begin{equation}\Pi_{02-}x(0)=c_{02-},\label{equ:25}\end{equation}
where $\Pi_{02+} x(\tau)$ is the vector with the first $p$ components of the vector $\Pi_{02}x(\tau)$, and $\Pi_{02-} x(\tau)$ is the vector containing other components of vector $\Pi_{02} x(\tau)$. Here $c_{02-}$ is a constant vector which will be found below.

It can be proved that $\Pi_{0}x(\tau)$ possesses the exponential estimate
\begin{equation}||\Pi_{0}x(\tau)||\leq \kappa \exp(-\alpha\tau),\;\tau\geq 0.\label{equ:26}\end{equation}

Here and below we will use $\kappa$ and $\alpha$ to represent appropriate positive numbers, which
are, generally speaking, different in different inequalities.

We can rewritten the system (\ref{equ:23}) for $\tau\geq\tau_0$ as follows
\begin{equation}\frac{d\Pi_{02} x}{d\tau}=\Lambda(0,0)\Pi_{02}x+l_{02}(\Pi_0x),\label{equ:27}\end{equation}
where $\Lambda(0,0)=diag\{\Lambda_+(0,0),\Lambda_-(0,0)\}$, $l_{02}(\Pi_{02}x)=f_2(\overline{x}_0(0)+\Pi_0x,0,0)-f'_{2x}(\overline{x}_0(0),0,0)\Pi_{02}x$, $f'_{2x}=\left(\frac{\partial f_{2i}}{\partial x_j}\right)_{i,j=\overline{2,n}}$. By construction $l_{02}(0)=0$.

Let us define vectors $l_{02+}(\Pi_{02}x)$ and $l_{02-}(\Pi_{02}x)$ in the same way as $\Pi_{02+} x$ and $\Pi_{02-} x$, respectively.

Instead of (\ref{equ:27}), (\ref{equ:24}), (\ref{equ:25}) we consider the equivalent system of integral equations
\begin{equation}\Pi_{02+} x(\tau)=-\int\limits_{\tau}^{\infty}\exp(\Lambda_+(0,0)(\tau-s))l_{02+}(\Pi_{02}x(s))ds,\label{equ:28}\end{equation}
\begin{equation}\Pi_{02-} x(\tau)=\exp(\Lambda_-(0,0)(\tau-\tau_0))\Pi_{02-} x(\tau_0)+\nonumber\end{equation}
\begin{equation}+
\int\limits_{\tau_0}^{\tau}\exp(\Lambda_-(0,0)(\tau-s))l_{02-}(\Pi_{02}x(s))ds.\label{equ:29}\end{equation}

For the proof of the estimate (\ref{equ:26}) we use the method of successive approximations to the system (\ref{equ:28}), (\ref{equ:29}). Let
\[\Pi_{02+}^{(0)} x(\tau)=0,\;\Pi_{02-}^{(0)} x(\tau)=\exp(\Lambda_-(0,0)(\tau-\tau_0))\Pi_{02-} x(\tau_0),\]
\[\Pi_{02+}^{(l)} x(\tau)=-\int\limits_{\tau}^{\infty}\exp(\Lambda_+(0,0)(\tau-s))l_{02+}(\Pi_{02}^{(l-1)}x(s))ds,\]
\[\Pi_{02-}^{(l)} x(\tau)=\exp(\Lambda_-(0,0)(\tau-\tau_0))\Pi_{02-} x(\tau_0)+\]
\[ +\int\limits_{\tau_0}^{\tau}\exp(\Lambda_-(0,0)(\tau-s))l_{02-}(\Pi_{02}^{(l-1)}x(s))ds,\;l\in N.\]

Let $\alpha_1$ and $r_1$ be positive numbers such that \[\exp(\Lambda_+(0,0)(\tau-s))\leq r_1\exp(\alpha_1(\tau-s)),\;\tau_0\leq\tau\leq s,\] \[\exp(\Lambda_-(0,0)(\tau-s))\leq r_1\exp(-\alpha_1(\tau-s)),\;\tau_0\leq s\leq\tau.\]
Note that $\alpha_1=\min\{\re \eta_1(0,0), -\re \eta_2(0,0)\}$.

Also, the conditions 6, 9 and 10 implies that for every $\delta>0$ there is a $\tau_0=\tau_0(\delta)$ such that \[||\Pi_0x(\tau_0)||\leq \delta\] for $\tau\geq \tau_0$.
Then for $\tau\geq \tau_0$ we have \[||\Pi_0^{(0)}x(\tau)||\leq r_1\exp(-\alpha_1(\tau-\tau_0))||\Pi_{02-}(\tau_0)||\leq r_1\delta \exp(-\alpha_1(\tau-\tau_0))\leq \] \[\leq r_1\delta\exp(-\alpha_0(\tau-\tau_0)),\;0<\alpha_0<\alpha_1.\]

According to the Lagrange finite-increments formula for every $\varepsilon>0$ there is a $\mu=\mu(\varepsilon)$ such that
\begin{equation}||l_{02}(u)-l_{02}(v)||\leq \varepsilon ||u-v||,\label{equ:30}\end{equation} for all $||u||\leq \mu$, $||v||\leq \mu$.

Put \[\delta_1=\varepsilon r_1\max\left\{\frac{1}{\alpha_1+\alpha_0},\frac{1}{\alpha_1-\alpha_0}\right\}.\]

We assume that a $\delta$ is so small that \[\frac{r_1\delta}{1-\delta_1}\leq \mu.\] Such a $\delta$ always exists for sufficiently large $\tau_0=\tau_0(\delta)$.
Then \[||\Pi_0^{(0)}x(\tau)||\leq \mu,\;\tau\geq\tau_0,\] and
\[||\Pi_{02+}^{(1)}x(\tau)-\Pi_{02+}^{(0)}x(\tau)||\leq \varepsilon r_1^2 \delta \exp(\alpha_1\tau+\alpha_0\tau_0)\int\limits_{\tau}^{\infty}\exp(-(\alpha_1 +\alpha_0)s)ds\leq\]
\[\leq r_1\delta\delta_1\exp(-\alpha_0(\tau-\tau_0)),\]
\[||\Pi_{02-}^{(1)}x(\tau)-\Pi_{02-}^{(0)}x(\tau)||\leq \varepsilon r_1^2 \delta \exp(-\alpha_1\tau+\alpha_0\tau_0)\int\limits_{\tau_0}^{\tau}\exp((\alpha_1 -\alpha_0)s)ds\leq\]
\[\leq r_1\delta\delta_1\exp(-\alpha_0(\tau-\tau_0)).\]
Thus, we hawe
\[||\Pi_{02}^{(1)}x(\tau)||\leq r_1\delta(1+\delta_1)\exp(-\alpha_0(\tau-\tau_0))\leq \frac{r_1\delta}{1-\delta_1}\exp(-\alpha_0(\tau-\tau_0))\leq \mu, \;\tau\geq \tau_0.\]

We prove, by induction, that
\[||\Pi_{02}^{(l)}x(\tau)-\Pi_{02}^{(l-1)}x(\tau)||\leq r_1\delta\delta_1^l\exp(-\alpha_0(\tau-\tau_0)),\;l\in N,\] and
\[||\Pi_{02}^{(l)}x(\tau)||\leq r_1\delta(1+\delta_1+...+\delta_1^l)\exp(-\alpha_0(\tau-\tau_0)), \;l\in N\cup \{0\},\; \tau\geq \tau_0.\]
Therefore, \[\underset{l\rightarrow +\infty}{\lim}\Pi_{02}^{(l)} x(\tau)=\Pi_{02} x(\tau)\] and \[||\Pi_{02}x(\tau)||\leq \frac{r_1\delta}{1-\delta_1}\exp(-\alpha_0(\tau-\tau_0)), \; \tau\geq \tau_0.\]
\cite{har,vas1}.

For $0\leq \tau\leq\tau_0$ the solution $\Pi_{02}x(\tau)$ is bounded by some constant $r_2$
\[||\Pi_{02}x(\tau)||\leq r_2.\]
If we put \[r=\max\left\{r_2\exp(\alpha_0\tau_0),\frac{r_1\delta}{1-\delta_1}\exp(\alpha_0\tau_0)\right\},\] then
\[||\Pi_{02}x(\tau)||\leq r\exp(-\alpha_0\tau), \; \tau\geq 0.\]
Therefore, the condition 9 implies that the estimate (\ref{equ:26}) is valid for $0<\alpha<\alpha_0$.

Note that the solution $\Pi_{02} x(\tau)$ of the system (\ref{equ:23}) satisfying (\ref{equ:24}), (\ref{equ:25}) can be written in the form
\begin{equation}\Pi_{02+} x(\tau)=-\int\limits_{\tau}^{\infty}\exp(\Lambda_+(0,0)(\tau-s))l_{02+}(\Pi_{02}x(s))ds,\label{equ:31}\end{equation}
\begin{equation}\Pi_{02-} x(\tau)=\exp(\Lambda_-(0,0)\tau) c_{02-}+\nonumber\end{equation}
\begin{equation}+
\int\limits_{0}^{\tau}\exp(\Lambda_-(0,0)(\tau-s))l_{02-}(\Pi_{02}x(s))ds.\label{equ:32}\end{equation}
This solution depends on $c_{02-}$ as a parameter, i.e.  $\Pi_{0} x(\tau) =\Pi_{0} x(\tau,c_{02-})$.

Let us now consider the system (\ref{equ:21}), which we will write as
\begin{equation}A(T,0)\frac{dQ_0x}{d\xi}=f'_x(\overline{x}_0(T),T,0)Q_0x+m_0(Q_0x),\label{equ:33}\end{equation}
where $m_0(Q_0x)=f(\overline{x}_0(T)+Q_0x,T,0)-f'_x(\overline{x}_0(T),T,0)Q_0x$.

It follows from condition 7 that there exists the nonsingular matrix $U$ such that
\[U^{-1}A^{-1}(T,0)f_x'(\overline{x}_0(T),T,0)U=W(T,0)\equiv\] \[\equiv diag\{w_1(T,0),w_2(T,0),...,w_n(T,0)\}.\]

Putting $Q_0x=UR_0x$ in (\ref{equ:33}), we rewrite the system (\ref{equ:33}) as
\begin{equation}\frac{dR_0x}{d\xi}=W(T,0)R_0x+p_0(R_0x),\label{equ:34}\end{equation}
where $p_0(R_0x)=U^{-1}A^{-1}(T,0)m_0(UR_0x)$.

Note that there is an inequality similar to (\ref{equ:30}) for the function $p_0(R_0x)$, i.e.,
for every $\varepsilon>0$ there exists a $\mu(\varepsilon)$ such that
\begin{equation}||p_0(u)-p_0(v)||\leq \varepsilon ||u-v||,\label{equ:35}\end{equation} for all $||u||\leq \mu$, $||v||\leq \mu$.

Let us denote by $R_{0+} x$ and $p_{0+}(R_0 x)$ the vectors with the first $p+1$ components of the vectors $R_0x$ and $p_0(R_0x)$, respectively, and through $R_{0-} x$ and $p_{0-}(R_0x)$ we denote the vectors containing other components of vectors $R_0 x$ and $p_0(R_0x)$.

Then the system (\ref{equ:34}) can be written in the form
\begin{equation}\frac{dR_{0+}x}{d\xi}=W_+(T,0)R_{0+}x+p_{0+}(R_0x),\label{equ:36}\end{equation}
\begin{equation}\frac{dR_{0-}x}{d\xi}=W_-(T,0)R_{0-}x+p_{0-}(R_0x),\label{equ:37}\end{equation}
where \[W_+(T,0)=diag\{w_1(T,0),w_2(T,0),...,w_{p+1}(T,0)\},\] \[W_-(T,0)=diag\{w_{p+2}(T,0),w_{p+2}(T,0),...,w_n(T,0)\}.\]
11. Assume that the system (\ref{equ:34}) has the solution $R_0x=R_0x(\xi)$ such that
\begin{equation}R_0x(\xi)\rightarrow 0,\;\xi\rightarrow -\infty,\label{equ:38}\end{equation}
and \begin{equation}R_{0+}x(0)=c_{0+},\label{equ:39}\end{equation}
where $c_{0+}$ is a constant vector which will be found below.

As before, we can prove the exponential estimates
\[||R_0 x(\xi)||\leq
\kappa\exp(\beta\xi),\;\xi\leq 0,\;0<\beta<\beta_1,\]
\[\exp(W_+(T,0)\xi)\leq \exp(\beta_1\xi),\;\exp(-W_-(T,0)\xi)\leq \exp(-\beta_1\xi),\;\xi\leq 0.\]

Note that this solution satisfies the system of the integral equations
\begin{equation}R_{0+} x(\xi)=\exp(W_+(T,0)\xi)c_{0+}-\int\limits_{\xi}^0\exp(W_+(T,0)(\xi-s))p_{0+}(R_0x)ds,\label{equ:40}\end{equation}
\begin{equation}R_{0-} x(\xi)=\int\limits_{-\infty}^{\xi}\exp(W_-(T,0)(\xi-s))p_{0-}(R_0x)ds.\label{equ:41}\end{equation}

The solution of (\ref{equ:34}) depends on $c_{0+}$ as a parameter, i.e.  $R_0 x(\xi) =R_0 x(\xi,c_{0+})$.

Let us now define the vectors $c_{02-}$ and $c_{0+}$. In a neighborhood of $t=0$ function $Qx(\xi,\varepsilon)$ is as small as, like a function $\Pi x(\tau,\varepsilon)$ in a neighborhood of $t=T$ \cite{vas1}. Therefore the equality (\ref{equ:16}) can be written in the form
\begin{equation}
  M(\overline{x}(0,\varepsilon)+\Pi x(0,\varepsilon))+N(\overline{x}(T,\varepsilon)+Q x(0,\varepsilon))=d(\varepsilon).
\nonumber\end{equation}
For $\varepsilon=0$ we get
\begin{equation}M(\overline{x}_0(0)+\Pi_0 x(0,c_{02-}))+N(\overline{x}_0(T)+UR_0x(0,c_{0+}))=d_0.\label{equ:42}\end{equation}
Assume that the following condition is satisfied.\\ 12. The equation (\ref{equ:42}) for $c_{02-}$, $c_{0+}$ have a solution.

{\bf Remark 1.} Let us find sufficient conditions for the existence $c_{02-}$ and $c_{0+}$.
Let \[M=\left(
        \begin{array}{ccc}
          M_{11} & M_{12} & M_{13} \\
          M_{21} & M_{22} & M_{23} \\
          M_{31} & M_{32} & M_{33} \\
        \end{array}
      \right),\;N=\left(
        \begin{array}{ccc}
          N_{11} & N_{12} & N_{13} \\
          N_{21} & N_{22} & N_{23} \\
          N_{31} & N_{32} & N_{33} \\
        \end{array}
      \right),\]
      \[K=\left(
        \begin{array}{ccc}
          K_{11} & K_{12} & K_{13} \\
          K_{21} & K_{22} & K_{23} \\
          K_{31} & K_{32} & K_{33} \\
        \end{array}
      \right),\;K=NU,\]
where the diagonal blocks $M_{11}$, $M_{22}$, and $M_{33}$ are rectangular matrices of dimensions $m_1\times 1$, $m_2\times p$, and $m_3\times q$, respectively. Here, of course, $m_1+m_2+m_3=n$. The matrices $N$ and $K$ have the same structure as $M$.

Then the equation (\ref{equ:42}) can be written as \[Dc_0=\varphi,\] where
\[c_0=\left(
        \begin{array}{c}
          c_{0+}\\
          c_{02-} \\
        \end{array}
      \right),\; D=\left(\begin{array}{ccc}
                            K_{11} & K_{12} & M_{13} \\
                            K_{21} & K_{22} & M_{23} \\
                            K_{31} & K_{32} & M_{33} \\
                            \end{array}\right),\]
\[\varphi=d_0-M\overline{x}_0(0)-N\overline{x}_0(T)-\left(
        \begin{array}{ccc}
          M_{11} & M_{12} & K_{13} \\
          M_{21} & M_{22} & K_{23} \\
          M_{31} & M_{32} & K_{33} \\
        \end{array}
      \right)\left(
                \begin{array}{c}
                  \Pi_{01}(\Pi_{02}x(0,c_{02-})) \\
                  \Pi_{02+}x(0,c_{02-}) \\
                  R_{0-}x(0,c_{0+}) \\
                \end{array}
              \right).\]

Suppose that the matrix $D$ is nonsingular. Then
\begin{equation}
  c_0=D^{-1}\varphi.
\label{equ:43}\end{equation}

According to the Lagrange finite-increments formula we get
\[||\Pi_0x(0,c_1)-\Pi_0x(0,c_2)||\leq \rho_1(c_1,c_2)||c_1-c_2||,\]
and
\[||R_0x(0,c_1)-R_0x(0,c_2)||\leq \rho_2(c_1,c_2)||c_1-c_2||,\;c_1,c_2\in R^n.\]
Then if \[||D^{-1}D_1||\left(2\underset{c_1,c_2\in R^n}{\sup}\, \rho_1(c_1,c_2)+\underset{c_1,c_2\in R^n}{\sup}\,\rho_2(c_1,c_2)\right)<1,\]
where \[D_1=\left(
        \begin{array}{ccc}
          M_{11} & M_{12} & K_{13} \\
          M_{21} & M_{22} & K_{23} \\
          M_{31} & M_{32} & K_{33} \\
        \end{array}
      \right),\]
the mapping (\ref{equ:43}) of the space $R^n$ into itself is a contraction mapping. Consequently, the equation (\ref{equ:43}) for $c_{02-}$, $c_{0+}$ has one and only one solution.

Equating coefficients of like powers of $\varepsilon$ in the equations (\ref{equ:13}) -- (\ref{equ:15}), we obtain
\begin{equation}f_x'(\overline{x}_0(t),t,0)\overline{x}_k=\sum\limits_{i=0}^{k-1}A_i(t)\frac{d\overline{x}_{k-i-1}(t)}{dt}-\overline{g}_k(t),\label{equ:44}\end{equation}
\begin{equation}A(0,0)\frac{d\Pi_kx}{d\tau}=f_x'(\overline{x}_0(0)+\Pi_0x(\tau),0,0)\Pi_kx(\tau)+r_k(\tau),\label{equ:45}\end{equation}
\begin{equation}A(T,0)\frac{dQ_kx}{d\xi}=f_x'(\overline{x}_0(T)+Q_0x(\xi),T,0)Q_kx(\xi)+q_k(\xi),\label{equ:46}\end{equation}
where
\[r_k(\tau)=g_k(\tau)-\sum\limits_{i=1}^k\Pi_iA(\tau)\frac{d\Pi_{k-i}x(\tau)}{d\tau},\;q_k(\xi)=h_k(\xi)-\sum\limits_{i=1}^kQ_iA(\xi)\frac{dQ_{k-i}x(\xi)}{d\xi}.\]

The conditions 6, 7 and 8 implies that $det f_x'(\overline{x}_0(t),t,0)\not=0$, $t\in[0;T]$. Then
\[\overline{x}_k(t)=(f_x'(\overline{x}_0(t),t,0))^{-1}\left(\sum\limits_{i=0}^{k-1}A_i(t)\frac{d\overline{x}_{k-i-1}(t)}{dt}-\overline{g}_k(t)\right),\;k\in N.\]

Let us consider the system (\ref{equ:45}). We set \[f_x'(\overline{x}_0(0)+\Pi_0x(\tau),0,0)=\left(
                                                   \begin{array}{cc}
                                                     C_1(\tau) & C_2(\tau) \\
                                                     C_3(\tau) & C_4(\tau) \\
                                                   \end{array}
                                                 \right),\]
where $C_4(\tau)$ is a square matrix of order $n-1$. Note that $f_x'(\overline{x}_0(0)+\Pi_0x(\tau),0,0)\rightarrow \Omega(0,0)$, $\tau\rightarrow +\infty$. Now the system (\ref{equ:45}) can be rewritten as
\begin{equation}\Pi_{k1}x=-\frac{1}{C_1(\tau)}\left(C_2(\tau)\Pi_{k2}x+r_{k1}(\tau)\right),\label{equ:47}\end{equation}
\begin{equation}\frac{d\Pi_{k2}x}{d\tau}=\left(C_4(\tau)-\frac{C_3(\tau)C_2(\tau)}{C_1(\tau)}\right)\Pi_{k2}x+r_{k2}(\tau)-
\frac{C_3(\tau)r_{k1}(\tau)}{C_1(\tau)},\label{equ:48}\end{equation}
where $\Pi_{k1}x$, $r_{k1}(\tau)$ are the first components of the vectors $\Pi_kx$, $r_k(\tau)$, and $\Pi_{k2}x$, $r_{k2}(\tau)$ are the rest components of the vectors $\Pi_kx$, $r_k(\tau)$. Note that $r_k(\tau)$ independent of $\Pi_kx$.

Let $\Phi(\tau)$ be a fundamental matrix solution of the corresponding homogeneous system
\begin{equation}\frac{d\Pi_{k2}x}{d\tau}=\Lambda(0,0)\Pi_{k2}x+\left(C_4(\tau)-\Lambda(0,0)-\frac{C_3(\tau)C_2(\tau)}{C_1(\tau)}\right)\Pi_{k2}x,\;\tau\geq \tau_0,\label{equ:49}\end{equation}
where $\Lambda(0,0)=diag\{\Lambda_+(0,0),\Lambda_-(0,0)\}$.

 Set \[\Phi(\tau)=\left(
                                                                                                            \begin{array}{cc}
                                                                                                             \Phi_{1}(\tau) & \Phi_{2}(\tau) \\
                                                                                                              \Phi_{3}(\tau) & \Phi_{4}(\tau) \\
                                                                                                            \end{array}
                                                                                                          \right),\]
where $\Phi_{1}(\tau)$ is the square matrix of order $p$.

Let $\varphi_{i}(\tau)$, $i=\overline{1,n-1}$, be the columns of the matrices $\Phi(\tau)$.
Then we have \cite{rap} \[\{\varphi_{i}(\tau)\}_j=\sigma_{ji}(\tau)\tau^{i-j}e^{\eta_1(0,0)(\tau-\tau_0)},\;i=\overline{1,p},\;j=\overline{1,n-1},\] and
\[\{\varphi_{i}(\tau)\}_j=\sigma_{ji}(\tau)\tau^{i-j}e^{\eta_2(0,0)(\tau-\tau_0)},\;i=\overline{p+1,n-1},\;j=\overline{1,n-1},\]
where
\[\sigma_{ji}(\infty)=\frac{1}{(i-j)!},\;j<i;\; i,j=\overline{1,p},\; \mbox{or}\; i,j=\overline{p+1,n-1},\]
\[\sigma_{jj}(\infty)=1,\; j=\overline{1,n-1},\]
\[\Phi_2(\infty)=0,\;\Phi_3(\infty)=0.\]
Here $\{\varphi_{i}(\tau)\}_j$ is the $j$th component of $\varphi_{i}(\tau)$.

Let \[\sigma_{ji}(\tau)=\frac{1}{(i-j)!}+\widetilde{\sigma}_{ji}(\tau),\;j<i,\;i,j=\overline{1,p},\; \mbox{or}\; i,j=\overline{p+1,n-1},\]
\[\sigma_{jj}(\tau)=1+\widetilde{\sigma}_{jj}(\tau),\; j=\overline{1,n-1}.\]

We set \[\Phi(\tau)=\left(\begin{array}{cc}
                               \Upsilon_{1}(\tau) & \Upsilon_{2}(\tau) \\
                               \Upsilon_{3}(\tau) & \Upsilon_{4}(\tau) \\
                                 \end{array}
                                 \right)\left(\begin{array}{cc}
                                \exp(\re\eta_1(0,0)\tau) I_p & 0 \\
                               0 & \exp(\re\eta_2(0,0)\tau) I_q \\
                                 \end{array}
                                 \right).\]
It follows from the method of construction of $\Phi(\tau)$ \cite{rap} that
\[|\widetilde{\sigma}_{ji}(\tau)|\leq \kappa\exp(-\alpha\tau),\;j\geq i,\;|\widetilde{\sigma}_{ji}(\tau)|\leq \frac{\kappa}{\tau},\;j< i,\]
$i,j=\overline{1,p}$, or $i,j=\overline{p+1,n-1}$, and
\[||\Upsilon_{2}(\tau)||\leq \kappa\exp(-\alpha\tau),\;||\Upsilon_{3}(\tau)||\leq \kappa\exp(-\alpha\tau),\;\tau\geq\tau_0.\]
Then the inverse matrix $\Phi^{-1}(\tau)$ can be write as
\[\Phi^{-1}(\tau)=\left(\begin{array}{cc}
                                \exp(-\re\eta_1(0,0)\tau) \Upsilon^{-1}_{1}(\tau) & 0 \\
                               0 & \exp(-\re\eta_2(0,0)\tau) \Upsilon^{-1}_{4}(\tau) \\
                                 \end{array}
                                 \right)\times\] \[\times(I_{n-1}+\Delta(\tau)),\]
where $||\Delta(\tau)||\leq \kappa\exp(-\alpha\tau)$, $\tau\geq \tau_0$.

Let $\widetilde{\Phi}(\tau)$, $\tau\geq 0$ be a fundamental matrix solution of the system (\ref{equ:49}) such that $\widetilde{\Phi}(0)=I_{n-1}$. Then for $\tau\geq \tau_0$ we have \[\widetilde{\Phi}(\tau)=\Phi(\tau)\Gamma_1,\] where $\Gamma_1$ is a constant square matrix.

The solution of the system (\ref{equ:48}) can be written in the following form
\[\Pi_{k2}x(\tau)=\widetilde{\Phi}(\tau)\left(\Gamma_1^{-1} a_k+\int\limits_0^{\tau}\widetilde{\Phi}^{-1}(k)l_{k2}(s)ds\right),\]
where \[l_{k2}(\tau)=r_{k2}(\tau)-
\frac{C_3(\tau)r_{k1}(\tau)}{C_1(\tau)},\] and $a_k$ is a $(n-1)$-dimensional constant vector which will be defined below.

If we set $\widetilde{\widetilde{\Phi}}(\tau)=\widetilde{\Phi}(\tau)\Gamma_1^{-1}$, the solution of the system (\ref{equ:48})
satisfying the conditions $\Pi_{k2+} x(\tau)\rightarrow 0$, $\tau\rightarrow +\infty$, and $\Pi_{k2-} x(0)=c_{k2-}$ can be written as
\[\Pi_{k2+}x(\tau)=-\widetilde{\widetilde{\Phi}}_{1}(\tau)\int\limits_{\tau}^{\infty}v_{k+}(s)ds+\widetilde{\widetilde{\Phi}}_{2}(\tau)\left(a_{k-}+\int\limits_0^{\tau}v_{k-}(s)ds\right),\]
\[\Pi_{k2-}x(\tau)=-\widetilde{\widetilde{\Phi}}_{3}(\tau)\int\limits_{\tau}^{\infty}v_{k+}(s)ds+\widetilde{\widetilde{\Phi}}_{4}(\tau)\left(a_{k-}+\int\limits_0^{\tau}v_{k-}(s)ds\right),\]
where
\[\widetilde{\widetilde{\Phi}}(\tau)=\left(\begin{array}{cc}
                               \widetilde{\widetilde{\Phi}}_{1}(\tau) & \widetilde{\widetilde{\Phi}}_{2}(\tau) \\
                               \widetilde{\widetilde{\Phi}}_{3}(\tau) & \widetilde{\widetilde{\Phi}}_{4}(\tau) \\
                                 \end{array}
                                 \right),\;a_k=\left(\begin{array}{c}
                               a_{k+} \\
                               a_{k-} \\
                                 \end{array}
                                 \right),\] \[v_k(\tau)= (\widetilde{\widetilde{\Phi}}(\tau))^{-1}l_{k2}(\tau),\;v_k(\tau)=\left(\begin{array}{c}
                               v_{k+}(\tau) \\
                               v_{k-}(\tau) \\
                                 \end{array}
                                 \right),\]
the vectors $\Pi_{k2+} x(\tau)$, $a_{k+}$, $v_{k+}(\tau)$ and $\Pi_{k2-} x(\tau)$, $a_{k-}$, $v_{k-}(\tau)$ have the same structure as $\Pi_{02+} x(\tau)$ and $\Pi_{02-} x(\tau)$, respectively. Here $\widetilde{\widetilde{\Phi}}(\tau)=\Phi(\tau)$ for $\tau\geq \tau_0$, and \[c_{k2-}=-\Gamma_{13}^{(-1)}\int\limits_0^{\infty}v_{k+}(s)ds+\Gamma_{14}^{(-1)}a_{k-},\] where
\[\Gamma_1^{-1}=\left(\begin{array}{cc}
                               \Gamma_{11}^{(-1)} & \Gamma_{12}^{(-1)} \\
                               \Gamma_{13}^{(-1)} & \Gamma_{14}^{(-1)} \\
                                 \end{array}
                                 \right).\]

Then, by construction, we obtain the exponential estimate
\[||\Pi_kx(\tau)||\leq \kappa \exp(-\alpha\tau),\;\tau\geq 0.\]

As before, the solution of (\ref{equ:47}), (\ref{equ:48}) depends on $a_{k-}$ as a parameter, i.e.  $\Pi_k x(\tau) =\Pi_k x(\tau,a_{k-})$.

Consider now the system (\ref{equ:46}). Let $Q_kx=UR_kx$. Then we can write the system (\ref{equ:46}) as
\begin{equation}\frac{dR_kx}{d\xi}=F(\xi)R_kx+p_k(\xi),\label{equ:50}\end{equation}
where \[F(\xi)=U^{-1}A^{-1}(T,0)f'_x(\overline{x}_0(T)+Q_0x(\xi),T,0)U,\] \[p_k(\xi)=U^{-1}A^{-1}(T,0)q_k(\xi).\]
Let \[F(\xi)=\left(
                                                   \begin{array}{cc}
                                                     F_1(\xi) & F_2(\xi) \\
                                                     F_3(\xi) & F_4(\xi) \\
                                                   \end{array}
                                                 \right),\]
where $F_1(\xi)$ is the square matrix of order $p+1$. Note that $F(\xi)\rightarrow W(T,0)$, $\xi\rightarrow -\infty$.

Let $\Psi(\xi)$ be a fundamental matrix solution of the corresponding homogeneous system
\begin{equation}\frac{dR_kx}{d\xi}=W(T,0)R_kx+(F(\xi)-W(T,0))R_kx(\xi),\;\xi\leq \xi_0.\label{equ:51}\end{equation}

We set \[\Psi(\xi)=\left(
                                                                                                            \begin{array}{cc}
                                                                                                             \Psi_{1}(\xi) & \Psi_{2}(\xi) \\
                                                                                                              \Psi_{3}(\xi) & \Psi_{4}(\xi) \\
                                                                                                            \end{array}
                                                                                                          \right),\]
where $\Psi_{1}(\xi)$ is the square matrix of order $p+1$.

Let $\psi_{i}(\xi)$, $i=\overline{1,n}$, be the columns of the matrix $\Psi(\xi)$.
Then we have \cite{rap}
\[\{\psi_{i}(\xi)\}_j=\sigma_{ji}(\xi)e^{w_i(T,0)(\xi-\xi_0)},\;i,j=\overline{1,n},\] where
\[\sigma_{ji}(-\infty)=0,\;i\not=j,\; \sigma_{jj}(-\infty)=1.\]
Here $\{\psi_{i}(\xi)\}_j$ is the $j$th component of $\psi_{i}(\xi)$.

Let \[\sigma_{jj}(\xi)=1+\widetilde{\sigma}_{jj}(\xi).\]
Then, by construction,
\[|\sigma_{ji}(\xi)|\leq \kappa\exp(\beta\xi),\;i\not =j, \;|\widetilde{\sigma}_{jj}(\xi)|\leq \kappa\exp(\beta\xi).\]

If we set \[\Psi(\xi)=\left(\begin{array}{cc}
                               \Theta_{1}(\xi) & \Theta_{2}(\xi) \\
                               \Theta_{3}(\xi) & \Theta_{4}(\xi) \\
                                 \end{array}
                                 \right)\exp(W(T,0)(\xi-\xi_0)),\]
the inverse matrix $\Psi^{-1}(\xi)$ can be write as
\[\Psi^{-1}(\xi)=\exp(-W(T,0)(\xi-\xi_0))(I_n+\Delta(\xi)),\]
where $||\Delta(\xi)||\leq \kappa\exp(\beta\xi)$, $\xi\leq \xi_0$.

Let $\widetilde{\Psi}(\xi)$, $\xi\leq 0$ be a fundamental matrix solution of the system (\ref{equ:51}) such that $\widetilde{\Psi}(0)=I_{n}$. Then, as before, for $\xi\leq \xi_0$ we have \[\widetilde{\Psi}(\xi)=\Psi(\xi)\Gamma_2,\] where $\Gamma_2$ is a constant square matrix.

The solution of the system (\ref{equ:50}) can be written in the form
\[R_{k}x(\xi)=\widetilde{\Psi}(\xi)\left(\Gamma_2^{-1} b_k+\int\limits_0^{\xi}\widetilde{\Psi}^{-1}(s)p_{k}(s)ds\right),\]
where $b_k$ is a $n$-dimensional constant vector which will be defined below.

If we set $\widetilde{\widetilde{\Psi}}(\xi)=\widetilde{\Psi}(\xi)\Gamma_2^{-1}$, the solution of the system (\ref{equ:50})
satisfying the conditions $R_{k+} x(0)=c_{k+}$, and $R_{k-} x(\xi)\rightarrow 0$, $\xi\rightarrow -\infty$, can be written as
\[R_{k+}x(\xi)=\widetilde{\widetilde{\Psi}}_{1}(\xi)\left(b_{k+}+\int\limits_0^{\xi}u_{k+}(s)ds\right)+\widetilde{\widetilde{\Psi}}_{2}(\xi)\int\limits_{-\infty}^{\xi}u_{k-}(s)ds,\]

\[R_{k-}x(\xi)=\widetilde{\widetilde{\Psi}}_{3}(\xi)\left(b_{k+}+\int\limits_0^{\xi}u_{k+}(s)ds\right)+\widetilde{\widetilde{\Psi}}_{4}(\xi)\int\limits_{-\infty}^{\xi}u_{k-}(s)ds,\]
where
\[\widetilde{\widetilde{\Psi}}(\xi)=\left(\begin{array}{cc}
                               \widetilde{\widetilde{\Psi}}_{1}(\xi) & \widetilde{\widetilde{\Psi}}_{2}(\xi) \\
                               \widetilde{\widetilde{\Psi}}_{3}(\xi) & \widetilde{\widetilde{\Psi}}_{4}(\xi) \\
                                 \end{array}
                                 \right),\;b_k=\left(\begin{array}{c}
                               b_{k+} \\
                               b_{k-} \\
                                 \end{array}
                                 \right),\] \[u_k(\xi)= (\widetilde{\widetilde{\Psi}}(\xi))^{-1}p_{k}(\xi),\;u_k(\xi)=\left(\begin{array}{c}
                               u_{k+}(\xi) \\
                               u_{k-}(\xi) \\
                                 \end{array}
                                 \right),\]
the vectors $R_{k+} x(\xi)$, $b_{k+}$, $u_{k+}(\xi)$ and $R_{k-} x(\xi)$, $b_{k-}$, $u_{k-}(\xi)$ have the same structure as $R_{0+} x(\xi)$ and $R_{0-} x(\xi)$, respectively. As before, $\widetilde{\widetilde{\Psi}}(\xi)=\Psi(\xi)$ for $\xi\leq \xi_0$, and \[c_{k+}=\Gamma_{21}^{(-1)}b_{k+}+\Gamma_{22}^{(-1)}\int\limits_{-\infty}^0u_{k-}(s)ds,\] where
\[\Gamma_2^{-1}=\left(\begin{array}{cc}
                               \Gamma_{21}^{(-1)} & \Gamma_{22}^{(-1)} \\
                               \Gamma_{23}^{(-1)} & \Gamma_{24}^{(-1)} \\
                                 \end{array}
                                 \right).\]

Then, by construction, we obtain the exponential estimate
\[||R_kx(\xi)||\leq \kappa \exp(\beta\xi),\;\xi\leq 0.\]

As before, the solution of (\ref{equ:50}) depends on $b_{k+}$ as a parameter, i.e.  $R_k x(\xi) =R_k x(\xi,b_{k+})$.

The system \begin{equation}M(\overline{x}_k(0)+\Pi_k x(0,c_{k2-}))+N(\overline{x}_k(T)+UR_kx(0,c_{k+}))=d_k\label{equ:52}\end{equation}
with respect to $c_{k2-}$ and $c_{k+}$ is structurally the same as system (\ref{equ:42}).

Assume that the following condition is satisfied.\\ 13. The system (\ref{equ:52}) for $c_{k2-}$, $c_{k+}$ have a solution.

\section{Asymptotic behaviour of the constructed formal solution}

Let's prove that the constructed formal solution of the problem (\ref{equ:10}), (\ref{equ:11}) has the asymptotic properties. For this purpose we make the substitution \begin{equation}x(t,\varepsilon)=x_l(t,\varepsilon)+y(t,\varepsilon),\label{equ:53}\end{equation}
in system (\ref{equ:10}), where
$x_l(t,\varepsilon)=\sum\limits_{k=0}^l\varepsilon^s(\overline{x}_k(t)+\Pi_kx(\tau)+Q_kx(\xi))$, and $y(t,\varepsilon)$ is a new unknown function. Then it follows from conditions 6 and 8 system (\ref{equ:10}) can be written in the form
\begin{equation}\varepsilon B(t,\varepsilon)\frac{dy}{dt}=y+g(y,t,\varepsilon),\label{equ:54}\end{equation}
where \[B(t,\varepsilon)=(f_x'(\overline{x}_0(t),t,0))^{-1}A(t,\varepsilon),\] \[\displaystyle g(y,t,\varepsilon)=(f_x'(\overline{x}_0(t),t,0))^{-1}\biggl(f(x_{l}(t,\varepsilon)+y,t,\varepsilon)-  \] \[-f_x'(\overline{x}_0(t),t,0)y- \varepsilon A(t,\varepsilon)\frac{dx_{l}(t,\varepsilon)}{dt}\biggr).\]

The boundary conditions for the system (\ref{equ:54}) take the form:
\begin{equation}My(0,\varepsilon)+Ny(T,\varepsilon)=\varepsilon^{l+1}n(\varepsilon),\label{equ:55}\end{equation} where $n(\varepsilon)=O(1)$, as $\varepsilon\rightarrow 0+$.

We put
\[\overline{x}_l(t,\varepsilon)=\sum\limits_{k=0}^l\varepsilon^k\overline{x}_k(t),\; \Pi_lx(\tau,\varepsilon)=\sum\limits_{k=0}^l\varepsilon^k\Pi_k(\tau),\; Q_l(\xi,\varepsilon)=\sum\limits_{k=0}^l\varepsilon^kQ_k(\xi).\]
Then
\[f(\overline{x}_l(t,\varepsilon)+\Pi_lx(\tau,\varepsilon)+Q_lx(\xi,\varepsilon),t,\varepsilon)=f(\overline{x}_l(t,\varepsilon),t,\varepsilon)+\]
\[+(f(\overline{x}_l(t,\varepsilon)+\Pi_lx(\tau,\varepsilon),t,\varepsilon)-f(\overline{x}_l(t,\varepsilon),t,\varepsilon))+\] \[+(f(\overline{x}_l(t,\varepsilon)+Q_lx(\xi,\varepsilon),t,\varepsilon)-f(\overline{x}_l(t,\varepsilon),t,\varepsilon))+\] \[+f(\overline{x}_l(t,\varepsilon)+\Pi_lx(\tau,\varepsilon)+Q_lx(\xi,\varepsilon),t,\varepsilon)-
f(\overline{x}_l(t,\varepsilon)+\Pi_lx(\tau,\varepsilon),t,\varepsilon)-\] \[ -f(\overline{x}_l(t,\varepsilon)+Q_lx(\xi,\varepsilon),t,\varepsilon)+f(\overline{x}_l(t,\varepsilon),t,\varepsilon).\]
Considering this expression on the segment $[0;\frac{T}{2}]$ and $[\frac{T}{2};T]$, we obtain
\[f(\overline{x}_l(t,\varepsilon)+\Pi_lx(\tau,\varepsilon)+Q_lx(\xi,\varepsilon),t,\varepsilon)=f(\overline{x}_l(t,\varepsilon),t,\varepsilon)+\] \[+(f(\overline{x}_l(t,\varepsilon)+\Pi_lx(\tau,\varepsilon),t,\varepsilon)-f(\overline{x}_l(t,\varepsilon),t,\varepsilon))+\] \[+(f(\overline{x}_l(t,\varepsilon)+Q_lx(\xi,\varepsilon),t,\varepsilon)-f(\overline{x}_l(t,\varepsilon),t,\varepsilon))+O(\varepsilon^{l+1})=\]
\[=\sum\limits_{k=0}^l\varepsilon^k(\overline{f}_k(t)+\Pi_kf(\tau)+Q_kf(\xi))+O(\varepsilon^{l+1}).\]

It follows from conditions 5 -- 7 that the sequence
\[e(w,t,\varepsilon), 1, ..., 1\] of $n$ polynomials is the set of all invariant factors of $B(t,\varepsilon)$ for each $t\in[0;T]$.
It is clear that $e(w,t,\varepsilon)$ is a polynomial of degree $n$. Therefore, the degrees of invariant factors of $B(t,\varepsilon)$ are constant on the segment $[0;T]$.

Let $U(t,\varepsilon)$ be a square matrix of orders $n$ such that
\[U(t,\varepsilon)= diag \{U_1(t,\varepsilon),U_p(t,\varepsilon),U_q(t,\varepsilon)\},\] where
 \[U_1(t,\varepsilon)=\theta_1(t,\varepsilon),\;U_p(t,\varepsilon)=\Theta_p(t,\varepsilon)+ N_p,\; U_q(t,\varepsilon)=\Theta_q(t,\varepsilon)+ N_q,\] \[\Theta_p(t,\varepsilon)=diag\{\theta_2(t,\varepsilon),\theta_3(t,\varepsilon),...,\theta_{p+1}(t,\varepsilon)\},\]
\[\Theta_q(t,\varepsilon)=diag\{\theta_{p+2}(t,\varepsilon),\theta_{p+3}(t,\varepsilon),...,\theta_n(t,\varepsilon)\},\]
and
\[\theta_1(0,0)=0,\]
\[\theta_i(t,0)=\frac{1}{w_i(t,0)}\rightarrow \frac{1}{\nu_1(0,0)},\;t\rightarrow 0,\;i=\overline{2,p+1},\]
\[\theta_i(t,0)=\frac{1}{w_i(t,0)}\rightarrow \frac{1}{\nu_2(0,0)},\;t\rightarrow 0,\;i=\overline{p+2,n}.\]

The matrices $B(t,\varepsilon)$ and $U(t,\varepsilon)$ are point-wise similar on the segment $[0;T]$.

Then there exists a nonsingular sufficiently smooth matrix $V_1(t,\varepsilon)$, $(t,\varepsilon)\in [0;t_1]\times [0;\varepsilon_1]$, $t_1\leq T$, $\varepsilon_1\leq \varepsilon_0$, \cite{was1} such that
  \[V_1^{-1}(t,\varepsilon)B(t,\varepsilon)V_1(t,\varepsilon)=U(t,\varepsilon).\]

The eigenvalues of $B(t,\varepsilon)$ are distinct for $t_0\leq t\leq T$, $t_0<t_1$. Thus there exists a nonsingular sufficiently smooth matrix $V_2(t,\varepsilon)$, $(t,\varepsilon)\in [t_0;T]\times [0;\varepsilon_1]$, \cite{gan} such that
\[V_2^{-1}(t,\varepsilon)B(t,\varepsilon)V_2(t,\varepsilon)=U(t,\varepsilon).\]

Define the matrix $V(t,\varepsilon)$ by \[V(t,\varepsilon)=e_1(t)V_1(t,\varepsilon)+e_2(t)V_2(t,\varepsilon),\] where
$e_1(t)$, $e_2(t)\in C^{\infty}[0;T]$, and \[e_1(t)\not=0,\;t\in [0;t_1),\;e_1(t)\equiv 0,\;t\in[t_1;T],\]
\[e_2(t)\equiv 0,\;t\in [0;t_0],\;e_2(t)\not=0,\;t\in(t_0;T].\]

We also assume that the elements of the matrices $V_1(t,\varepsilon)$ and $V_2(t,\varepsilon)$ are sufficiently smooth on the segment $[0;T]$.
Then $V(t,\varepsilon)$ is the nonsingular sufficiently smooth matrix for $0\leq t\leq T$. In fact,
\[V(t,\varepsilon)=e_1(t)V_1(t,\varepsilon),\;t\in[0;t_0],\]
\[V(t,\varepsilon)=e_2(t)V_2(t,\varepsilon),\;t\in[t_1;T].\]

Consider the matrix $V(t,\varepsilon)$ on the interval $(t_0;t_1)$. On this interval
\[det V(t,\varepsilon)=det(e_1(t)V_1(t,\varepsilon)+e_2(t)V_2(t,\varepsilon))=\] \[=e_1^n(t)v_n(t,\varepsilon)
+e_1^{n-1}(t)e_2(t)v_{n-1}(t,\varepsilon)+...+e_1(t)e_2^{n-1}(t)v_1(t,\varepsilon)+e_2^n(t)v_0(t,\varepsilon),\]
where $v_n(t,\varepsilon)=det V_1(t,\varepsilon)$, $v_0(t,\varepsilon)=det V_2(t,\varepsilon)$, $v_i(t,\varepsilon)\in C[t_0;t_1]$, $i=\overline{1,n-1}$.
Let $e(t)=\frac{e_1(t)}{e_2(t)}$. Then we have
\[det V(t,\varepsilon)=e_2^n(t)(e^n(t)v_n(t,\varepsilon)+e^{n-1}(t)v_{n-1}(t,\varepsilon)+...+e(t)v_1(t,\varepsilon)+v_0(t,\varepsilon)).\]
Note that \[e^n(t)v_n(t,\varepsilon)+e^{n-1}(t)v_{n-1}(t,\varepsilon)+...+e(t)v_1(t,\varepsilon)+v_0(t,\varepsilon)\not=0,\;t\in(t_0;t_1),\]
since $v_0(t,\varepsilon)\not= 0$, $t\in(t_0;t_1)$. Thus $det V(t,\varepsilon)\not=0$, $t\in(t_0;t_1)$, and therefore $det V(t,\varepsilon)\not=0$, $t\in[0;T]$.

Putting $y=V(t,\varepsilon)z$ in (\ref{equ:54}), we rewrite the system (\ref{equ:54}) in the form
\begin{equation}\varepsilon U(t,\varepsilon)\frac{dz}{dt}=z+h(z,t,\varepsilon),\label{equ:56}\end{equation}
where $\displaystyle h(z,t,\varepsilon)=V^{-1}(t,\varepsilon)g(V(t,\varepsilon)z,t,\varepsilon)- \varepsilon U(t,\varepsilon)V^{-1}(t,\varepsilon)V'(t,\varepsilon)z$.
Note that
\[||h(z_1,t,\varepsilon)-h(z_2,t,\varepsilon)||\leq k_1(\varepsilon+\exp(-\alpha t/\varepsilon)+\exp(\beta(t-T)/\varepsilon))||z_1-z_2||\]
and \[||h(0,t,\varepsilon)||\leq k_2\varepsilon^{l+1},\;t\in [0;T],\]
for all $z_1$, $z_2\in D_{l+1}$, $D_{l+1}=\{z(t,\varepsilon)\in C[0;T]:||z(t,\varepsilon)||\leq k\varepsilon^{l+1}\}$.\\
14. Assume that $\re \theta_1(t,\varepsilon)>0$, $(t,\varepsilon)\in (0;T]\times (0;\varepsilon_0]$.

We will seek a bounded solution of the system (\ref{equ:56}). Such a solution must satisfy the system of integral equations
\begin{equation}z_1(t,\varepsilon)=
\exp\left(\frac{1}{\varepsilon}\int\limits_T^t\frac{du}{\theta_1(u,\varepsilon)}\right)\omega_1(\varepsilon)\varepsilon^{l+1}-\nonumber\end{equation}
\begin{equation}-\frac{1}{\varepsilon}\int\limits_t^T\frac{1}{\theta_1(s,\varepsilon)}
\exp\left(\frac{1}{\varepsilon}
\int\limits_s^t\frac{du}{\theta_1(u,\varepsilon)}\right)h_1(z(s,\varepsilon),s,\varepsilon)ds,\label{equ:57}\end{equation}
\begin{equation}z_{2+}(t,\varepsilon)=Z_{2+}(t,\varepsilon)Z_{2+}^{-1}(T,\varepsilon)\omega_{2+}(\varepsilon)\varepsilon^{l+1}-\nonumber\end{equation}
\begin{equation}-\frac{1}{\varepsilon}
\int\limits_t^TZ_{2+}(t,\varepsilon)Z_{2+}^{-1}(s,\varepsilon)U^{-1}_p(s,\varepsilon)h_{2+}(z(s,\varepsilon),s,\varepsilon)ds,\label{equ:58}\end{equation}
\begin{equation}z_{2-}(t,\varepsilon)=Z_{2-}(t,\varepsilon)Z_{2-}^{-1}(0,\varepsilon)\omega_{2-}(\varepsilon)\varepsilon^{l+1}+\nonumber\end{equation}
\begin{equation}+\frac{1}{\varepsilon}\int\limits_0^tZ_{2-}(t,\varepsilon)Z_{2-}^{-1}(s,\varepsilon)
U^{-1}_q(s,\varepsilon)h_{2-}(z(s,\varepsilon),s,\varepsilon)ds,\label{equ:59}\end{equation}
where $z(t,\varepsilon)=col(z_1(t,\varepsilon),z_{2+}(t,\varepsilon),z_{2-}(t,\varepsilon))$, $\omega(\varepsilon)=col (\omega_1(\varepsilon),\omega_{2+}(\varepsilon),\omega_{2-}(\varepsilon))$ is a constant vector such that $||\omega(\varepsilon)||=O(1)$ as $\varepsilon\rightarrow 0+$, $Z_{2+}(t,\varepsilon)$ and
$Z_{2-}(t,\varepsilon)$ are fundamental matrices solutions of
the systems \[\varepsilon
\frac{dz_{2+}}{dt}=U_p^{-1}(t,\varepsilon)z_{2+},\]
and \[\varepsilon
\frac{dz_{2-}}{dt}=U_q^{-1}(t,\varepsilon)z_{2-},\]
respectively.

Note that there exist positive numbers $\gamma$, $\zeta_+$, and $\zeta_-$ \cite{vas1} such that
\[||Z_{2+}(t,\varepsilon)Z_{2+}^{-1}(s,\varepsilon)||\leq \gamma\exp\left(\frac{\zeta_+(t-s)}{\varepsilon}\right),\;0\leq t\leq s\leq T,\]
\[||Z_{2-}(t,\varepsilon)Z_{2-}^{-1}(s,\varepsilon)||\leq \gamma\exp\left(-\frac{\zeta_-(t-s)}{\varepsilon}\right),\;0\leq s\leq t\leq T.\]

Assume that the following condition is satisfied.\\ 15. \[2k_1<1,\] \[2\gamma k_1k_3\max\left\{\frac{1}{\zeta_+}, \frac{1}{\zeta_-}\right\}<1,\] where \[\underset{t\in[0;T]}{\max}\{||U^{-1}_p(t,0)||,||U^{-1}_q(t,0)||\}\leq k_3.\]

Then the operator given by (\ref{equ:54}) -- (\ref{equ:56}) maps the set $D_{l+1}$ into itself for fixed $\omega(\varepsilon)$. This mapping is a contraction mapping. Consequently, the system (\ref{equ:54}) -- (\ref{equ:56}) has one and only one solution on the set $D_{l+1}$ for sufficiently large $k_3$. That is why the system (\ref{equ:53}) has unique solution $z=z(t,\varepsilon)$ as well. Furthermore, $||z(t,\varepsilon)||\leq k\varepsilon^{l+1}$, $t\in[0;T]$. Note that $z(t,\varepsilon)=z(t,\omega(\varepsilon),\varepsilon)$.

Finally, we take a constant vector $\omega$ so that the equality (\ref{equ:55}) hold. For this we will assume that the following condition is satisfied.\\
16. $det \widetilde{D}(0) \not=0$, where
\[\widetilde{D}(\varepsilon)=\left(
        \begin{array}{ccc}
          \widetilde{N}_{11}(T,\varepsilon) & \widetilde{N}_{12}(T,\varepsilon) & \widetilde{M}_{13}(0,\varepsilon) \\
          \widetilde{N}_{21}(T,\varepsilon) & \widetilde{N}_{22}(T,\varepsilon) & \widetilde{M}_{23}(0,\varepsilon) \\
          \widetilde{N}_{31}(T,\varepsilon) & \widetilde{N}_{32}(T,\varepsilon) & \widetilde{M}_{33}(0,\varepsilon) \\
        \end{array}
      \right),\] \[\widetilde{M}(0,\varepsilon)=e_1(0)MV_1(0,\varepsilon),\;\widetilde{N}(T,\varepsilon)=e_2(T)NV_2(T,\varepsilon).\]
Here, as before, the diagonal blocks $\widetilde{N}_{11}(T,\varepsilon)$, $\widetilde{N}_{22}(T,\varepsilon)$ and $\widetilde{M}_{33}(0,\varepsilon)$ are rectangular matrices of dimension $m_1\times 1$, $m_2\times p$, and $m_3\times q$, respectively.

Then from (\ref{equ:55}) we obtain
\begin{equation}
  \omega(\varepsilon)=\widetilde{D}^{-1}(\varepsilon)\widetilde{\varphi}(\varepsilon)
\label{equ:60},\end{equation} where
\[\widetilde{\varphi}(\varepsilon)=n(\varepsilon)-\frac{1}{\varepsilon^{l+1}}\left(
        \begin{array}{ccc}
          \widetilde{M}_{11}(0,\varepsilon) & \widetilde{M}_{12}(0,\varepsilon) & \widetilde{N}_{13}(T,\varepsilon) \\
          \widetilde{M}_{21}(0,\varepsilon) & \widetilde{M}_{22}(0,\varepsilon) & \widetilde{N}_{23}(T,\varepsilon) \\
          \widetilde{M}_{31}(0,\varepsilon) & \widetilde{M}_{32}(0,\varepsilon) & \widetilde{N}_{33}(T,\varepsilon) \\
        \end{array}
      \right)\left(
                \begin{array}{c}
                  z_1(0,\omega(\varepsilon),\varepsilon) \\
                  z_{2+}(0,\omega(\varepsilon),\varepsilon) \\
                  z_{2-}(T,\omega(\varepsilon),\varepsilon) \\
                \end{array}
              \right).\]

According to the Lagrange finite-increments formula we get
\[||z(t,\omega_1,\varepsilon)-z(t,\omega_2,\varepsilon)||\leq k_4\varepsilon^{l+1}||\omega_1-\omega_2||,\;t\in [0;T],\]
for all sufficiently small $||\omega_1||$ and $||\omega_2||$.\\
17. Assume that
\[2||\widetilde{D}^{-1}(0)\widetilde{D}_1(0)||
k_1k_4\max\left\{1,k_3\gamma\max\left\{\frac{1}{\xi_+},\frac{1}{\xi_-}\right\}\right\}<1,\] where
\[\widetilde{D}_1(\varepsilon)=\left(
        \begin{array}{ccc}
          \widetilde{M}_{11}(0,\varepsilon) & \widetilde{M}_{12}(0,\varepsilon) & \widetilde{N}_{13}(T,\varepsilon) \\
          \widetilde{M}_{21}(0,\varepsilon) & \widetilde{M}_{22}(0,\varepsilon) & \widetilde{N}_{23}(T,\varepsilon) \\
          \widetilde{M}_{31}(0,\varepsilon) & \widetilde{M}_{32}(0,\varepsilon) & \widetilde{N}_{33}(T,\varepsilon) \\
        \end{array}
      \right).\]

Then the mapping given by (\ref{equ:60}) of the space $R^n$ is a contraction mapping. Consequently, the equation (\ref{equ:60}) for $\omega(\varepsilon)$ has one and only one solution.

Thus, the main result of the paper can be formulated as follows.

{\bf Theorem.} If $A(t,\varepsilon)\in C^{l+1}(G)$, $f(x,t,\varepsilon)\in C^{l+1}(K)$ and the assumptions 3 -- 17 (the assumption 13 for $k=\overline{1,l}$) are satisfied, then there exists a unique solution $x=x(t,\varepsilon)$ of the boundary-value problem (\ref{equ:10}), (\ref{equ:11}) for sufficiently small $\varepsilon$, $0<\varepsilon\leq \varepsilon_1\leq \varepsilon_0$, such that
 \[||x(t,\varepsilon)-x_l(t,\varepsilon)||=O(\varepsilon^{l+1}),\;t\in[0;T],\;\varepsilon\rightarrow 0+.\]


\end{document}